\documentclass[centertags,12pt]{article}

\usepackage{amsmath}
\usepackage{amsthm}
\usepackage{amscd}
\usepackage{amssymb}
\usepackage{verbatim}

\usepackage{epsfig}
\usepackage{psfig}

\newcommand{\Q}{\mathbb Q}
\newcommand{\Z}{\mathbb Z}

\newcommand{\U}{\mathcal U}
\newcommand{\sg}{\Sigma_g}

\renewcommand{\ker}{K}
\newcommand{\kg}{{\mathcal K}_g}
\renewcommand{\lg}{{\mathcal L}_g}

\newcommand{\lt}{{\mathcal L}}
\newcommand{\kgs}{{\mathcal K}_{g,*}}
\newcommand{\tgs}{{\T}_{g,*}}
\newcommand{\tg}{{\T}_g}

\theoremstyle{plain}
\def\endproof{$\diamond$ \medskip}

\newtheorem{theorem}{Theorem}[section]
\newtheorem{proposition}[theorem]{Proposition}

\newtheorem{corollary}[theorem]{Corollary}

\newtheorem{question}[theorem]{Question}

\theoremstyle{definition}
\newtheorem{definition}[theorem]{Definition}

\usepackage{amsmath, amsthm, amssymb}
\usepackage{amsfonts}
\usepackage{euscript}
\usepackage{latexsym}
\usepackage{amscd, xy}
\usepackage{psfig}
\usepackage{epsfig}
\xyoption{all}
\CompileMatrices

\DeclareMathOperator{\Mod}{Mod}

\DeclareMathOperator{\T}{{\mathcal I}}
\DeclareMathOperator{\K}{{\mathcal K}}
\DeclareMathOperator{\SL}{SL_2}
\DeclareMathOperator{\GL}{GL}
\DeclareMathOperator{\Aut}{Aut}

\title{$\kg$ is not finitely generated\\}
\author{Daniel Biss\thanks{This research was conducted during the
period the first author served as a Clay Mathematics Institute
Long-Term Prize Fellow.} \ and Benson Farb\thanks{Supported in part 
by NSF grants DMS-9704640 and DMS-0244542.}}

\renewcommand\to{\longrightarrow}

\begin{document}

\maketitle

\section{Introduction}
Let $\Sigma_g$ be a closed orientable surface of genus $g$.  The {\em
mapping class group\/} $\Mod_g$ of $\Sigma_g$ is defined to be the group
of isotopy classes of orientation-preserving diffeomorphisms
$\Sigma_g\rightarrow \Sigma_g$.  Recall that an essential simple closed
curve $\gamma$ in $\Sigma_g$ is called a {\em bounding curve}, or
{\em separating curve}, if it is null-homologous in $\Sigma_g$
or, equivalently, if $\gamma$ separates $\Sigma_g$ into two
connected components.

Let $\K_g$ denote the subgroup of $\Mod_g$ generated by the
(infinite) collection of Dehn twists about bounding curves in
$\Sigma_g$.  Note that $\K_1$ is trivial.  
It has been a long-standing problem in the
combinatorial topology of surfaces to determine whether or not
the group $\K_g$ is finitely generated for $g\geq 2$.  For a discussion 
of this problem, see, e.g., \cite{Jo1,Jo3,Bi,Mo1,Mo3,Ak}.

McCullough-Miller \cite{MM} proved that $\K_2$ is not
finitely generated; Mess then proved that $\K_2$ is in fact an
infinite rank free group.  Akita proved in \cite{Ak} that 
for all $g\geq 2$, the rational homology 
$H_\ast(\K_g;\Q)$ is infinite-dimensional as a vector space over
$\Q$.  Note that since $\K_g$ admits a free action on the Teichmuller
space of $\Sigma_g$, which is contractible and finite-dimensional,
$\K_g$ has finite cohomological dimension.

For some time it was not known
if $\K_g$ was equal to, or perhaps a finite index subgroup of, the {\em
Torelli group} $\T_g$, which is the subgroup of
elements of $\Mod_g$ which act trivially on $H_1(\Sigma_g;\Z)$.
Powell \cite{Po} proved that $\K_2=\T_2$.  
Johnson proved in \cite{Jo2} that for $g\geq 3,$ the group $\K_g$ 
has infinite index in $\T_g$; he did this by constructing what is now called
the \emph{Johnson homomorphism}, which is the quotient map in the short
exact sequence
$$1\to\kg\to\T_g\to\wedge^3H/H \to 1$$
where $H=H_1(\Sigma_g;\Z).$  Johnson then proved in \cite{Jo3} that
$\T_g$ is 
finitely generated for all $g\geq 3$.  Our main result is the following.

\begin{theorem}
\label{theorem:main}
The group $\K_g$ is not finitely generated for any $g\geq 2$.  
\end{theorem}

We will also prove along the way that 
that the once-punctured analogue of $\K_g$ is 
not finitely generated.

Theorem \ref{theorem:main} answers Problem 10 of \cite{Mo2}, Problem
2.2(i) of \cite{Mo3}, and the question/conjecture on page 24 of \cite{Bi}.
We would still, however, like to know the answer to the following
question, asked by Morita (see \cite{Mo3}, Problem 2.2(ii)).

\begin{question}
Is $H_1(\K_g;\Z)$ finitely generated for $g\geq 3$?
\end{question}

Note that Birman-Craggs-Johnson (see, e.g. \cite{BC,Jo1}) and 
Morita \cite{Mo4} have found large abelian quotients of
$\K_g$.  We would also like to remark that Morita has 
discovered (see, e.g., \cite{Mo4,Mo2,Mo3}) a strong
connection between the algebraic structure of $\K_g$ and the Casson
invariant for homology $3$-spheres.  For example, Morita proved in
\cite{Mo4} that every integral homology 3-sphere can be obtained by
gluing two handlebodies along their boundaries via a map in $\K_g$;
further, he has been able to express the Casson invariant as a
homomorphism $\K_g\to\Z$ (see, e.g., \cite{Mo1}).

\bigskip
\noindent
{\bf Rough outline of the proof. }Our proof owes a great intellectual debt to
the paper \cite{MM} by D. McCullough and A. Miller, where the theorem is
demonstrated in the genus $2$ case; indeed we follow the same outline as their
proof.  

First, we find an action of $\K_g$ on the first homology of an 
abelian cover $Y$ of $\Sigma_g$ with Galois group $\Z^{2g-2}$.  
While $H_1(Y;\Z)$ is
infinitely generated, it is finitely generated as a module over the
group-ring of the Galois group of the cover.  We view this group-ring as
the ring $\lg$ of integral Laurent series in $2g-2$
variables. This action a priori gives a rather complicated
high-dimensional representation of 
$\K_g$.  We first project to a Laurent series ring $\lt$ in just one
variable, and then are able to
 find and quotient out a codimension two fixed submodule.  
This reduction to a $2$-dimensional representation is 
crucial for what follows.  We then analyze this representation
$$\rho:\K_g\to \SL(\lt)$$

The ring $\lt$ comes equipped with a discrete valuation, and so
$\SL(\lt)$ can be realized via Bruhat-Tits theory as a group of
automorphisms of a certain simplicial tree.  The Bass-Serre theory of
graphs 
of groups---equivalently, of groups acting on trees---is especially suited to
understanding whether or not such a group is finitely generated; one such
criterion is proven in \cite{MM}.  To complete the proof, we compute 
enough about the image of 
$\rho$ to apply this criterion to show that $\K_g$ is not finitely
generated.

\bigskip
\noindent
{\bf Acknowledgments. }  We would like to thank T.\ Cochran, F.\
Grunewald and S.\
Morita for useful comments on an earlier manuscript.

\section{Representing $\kg$ on an abelian cover}

Consider a standard symplectic
basis $\{a_1,\dots,a_g,b_1,\dots,b_g\}$ for $H_1(\sg;\Z),$ where $a_i
\cdot b_j = \delta_{i,j}$ and $a_i\cdot a_j = b_i \cdot b_j = 0.$
Here and throughout this article, the symbol $\cdot$ is used to
denote the algebraic intersection number of simple closed curves (or
homology classes).  By abuse of notation, we will
also sometimes view the $a_i$ and $b_i$ as elements of
$\pi_1(\Sigma_g)$, considered as relative to a fixed basepoint.

\subsection{The abelian cover}

Consider the free abelian
group $\Z^{2g-2}$ 
with generators $\{s_2,\dots,s_g,t_2,\dots,t_g\}$ and the surjection
$\psi:H_1(\sg;\Z)\rightarrow\Z^{2g-2}$ defined by 
$$\begin{array}{l}
\psi(a_1) = \psi(b_1) = 0\\
\psi(a_i) = s_i, \ i\geq 2\\
\psi(b_i) = t_i, \ i\geq 2
\end{array}
$$
Composing with the Hurewicz map $\pi_1(\sg)\rightarrow
H_1(\sg;\Z)$ gives a surjection
$\varphi:\pi_1(\sg)\rightarrow\Z^{2g-2};$ we denote the kernel of
$\varphi$ by $\ker$.

Let $p:Y\rightarrow\sg$ denote the covering corresponding to the
subgroup $\ker\subset\pi_1(\sg).$  The group $\Z^{2g-2}$ then acts 
on $Y$ by deck transformations.  This action induces an action of
$\Z^{2g-2}$ on $H_1(Y;\Z)$, which is consequently a $\Z\left[s_2^{\pm
1},\dots,s_g^{\pm 1},t_2^{\pm 1},\dots,t_g^{\pm 1}\right]$-module.  We
denote this Laurent series ring by $\lg.$  

It is rather easy to construct the cover $Y$ explicitly.  To this end,
consider the decomposition of $\sg$ into two subsurfaces 
$\Sigma_{g-1,1}$ and $\Sigma_{1,1}$ of genus $g-1$ and $1$,
respectively, obtained by cutting along the bounding curve
representing the homotopy class $[a_1,b_1].$  Note that the subspace
$H_1(\Sigma_{1,1};\Z)\subset 
H_1(\Sigma_g;\Z)$ is the span of $\{a_1,b_1\}$.  Let 
$Y'$ denote the universal abelian cover of
$\Sigma_{g-1,1}$, that is, the cover corresponding to the commutator
subgroup of $\pi_1(\Sigma_{g-1,1})$.  Since the boundary of $\Sigma_{g-1,1}$ is
null-homologous, it lifts to a collection of simple closed boundary
curves in $Y'$, indexed by the set $\Z^{2g-2}.$ We then obtain $Y$ by gluing
$\Sigma_{1,1}$ to each of these curves along its boundary.

The $\lg$-module structure of $H_1(Y;\Z)$ can now be read off from this
geometric description of $Y$.

\begin{proposition}
\label{Pidenthom}
The homology group $H_1(Y;\Z)$ is generated as an $\lg$-module by the
following ${{2g-2}\choose{2}} +1$ elements:
$a_1, b_1, [a_i,a_j]$ and $[b_i,b_j]$ for $2\leq i< j\leq
g,$ and $[a_i,b_j]$ for $2\leq i,j\leq g,$ with $[a_g,b_g]$ excepted.
Denote by $W$ the submodule of $H_1(Y;\Z)$
obtained by omitting the generators $a_1$ and $b_1$ from this list.
We then have
$p_*(a_1) = a_1, p_*(b_1) = b_1,$ and $p_*(c) = 0$ for any
$c\in W.$  Moreover, $$\frac{H_1(Y;\Z)}{W}$$ is a free $\lg$-module on
$\{a_1,b_1\}.$ 
\end{proposition}

\begin{proof}
It is a standard fact that the homology of the
surface obtained by sewing in discs along the boundary circles of
$Y'$ is generated as an $\lg$-module by the elements
$[a_i,a_j],$ $[b_i,b_j],$ and $[a_i,b_j]$ (for the sake of
normalization, we choose a single connected fundamental domain $X$ for
the 
action of $K$ on $Y$ and demand that all these generators be supported
in $X$).   Note that the 
resulting space is just the universal abelian cover of
$\Sigma_{g-1}$. The element
$[a_g,b_g]$ is omitted because the relation
$[a_2,b_2]\cdot\dots\cdot[a_g,b_g]$ in $\pi_1(\Sigma_{g-1})$ implies
that, in the homology of the cover, $[a_g,b_g]$ is in the span of the
$[a_i,b_i]$ for $2\leq i\leq g-1.$  

The identification of the images of the generators under
$p_*$ follows directly from their definition.
Finally, to compute $H_1(Y;\Z)/W,$ notice that $W$ is the image of the
natural map $H_1(Y';\Z)\rightarrow H_1(Y;\Z)$ and that the map
$H_1(Y,Y';\Z)\rightarrow H_0(Y';\Z)$ is zero; then observe that 
$Y/Y'$ is a wedge of tori, one for each element of $\Z^{2g-2}.$
\end{proof}

We will need to compute the algebraic intersection numbers of
certain curves in $Y$.  To ease the exposition of the next result, it
will be convenient to introduce another piece of notation.  We denote
the set $\{a_2,\dots,a_g,b_2,\dots,b_g\}$ by $\{c_1,\dots,c_{2g-2}\}$
and the set 
$\{s_2,\dots,s_g,t_2,\dots,t_g\}$ by $\{u_1,\dots,u_{2g-2}\}.$ Thus,
$H_1(Y)$ is generated as an $\lg$-module by the elements $a_1,b_1,$ and
$[c_i,c_j]$ for $1\leq i<j\leq 2g-2$ except $i = g$ and $j = 2g.$

\begin{proposition}
Suppose $i,j,i',j'\in\{1,\dots, 2g-2\}$ with $i\neq j$ and $i'\neq
j'.$  Assume
first that $\{i,j\}\cap\{i',j'\} = \emptyset.$ Then there
exists $\epsilon_{i,j,i',j'}\in\{1,0,-1\}$ such that 
$$[c_i,c_j]\cdot \left(u_1^{r_1}\cdots
u_{2g-2}^{r_{2g-2}}[c_{i'},c_{j'}]\right) = \hspace{3.2in}$$
\begin{equation}\label{eq:calculation}
= \begin{cases}
(-1)^{r_{i}+r_j+r_{i'}+r_{j'}}\epsilon_{i,j,i',j'} & r_k
\in\{\delta_{i,k},\delta_{j,k},-\delta_{i',k},-\delta_{j',k}\} 
\mbox{ for all $k$}\\ 
0 & \mbox{otherwise.}
\end{cases}\end{equation}
Now, assume that
$i = i'.$  Then there exists
$\epsilon_{i,j,j'}\in\{1,0,-1\}$ such 
that
$$[c_i,c_j]\cdot \left(u_1^{r_1}\cdots
u_{2g-2}^{r_{2g-2}}[c_{i},c_{j'}]\right) = \hspace{3.2in}$$
\begin{equation}\label{eq:calculationparallel}
= \begin{cases}
(-1)^{r_{i}+r_j+r_{j'}}\epsilon_{i,j,j'} & r_k \in 
\{\delta_{i,k},\delta_{j,k},-\delta_{i,k},-\delta_{j',k}\} \mbox{ for
all $k$} \\
0 & \mbox{otherwise.}
\end{cases}\end{equation}
Lastly,
\begin{equation}\label{eq:vanishing}
[c_i,c_j]\cdot \left(u_1^{r_1}\cdots
u_{2g-2}^{r_{2g-2}} a_1\right) =
[c_i,c_j]\cdot \left(u_1^{r_1}\cdots
u_{2g-2}^{r_{2g-2}} b_1\right) = 0
\end{equation}
regardless of the integers $r_k.$
\end{proposition}

\begin{proof}
Equation~(\ref{eq:vanishing}) is clear since the curves in question are 
disjoint.  To prove equation~(\ref{eq:calculation}), notice that
the curve representing the cycle $[c_i,c_j]$ is a kind of quadrilateral
beginning at some basepoint $y$ in the fundamental domain $X,$ then
passing to $u_i y,$ followed by $u_i u_j y$, then $u_j y,$ and then back
to the original basepoint $y$.  The curve $[c_{i'},c_{j'}]$ thus
intersects $[c_i,c_j]$ only once, at $y$, but this intersection is not
necessarily transverse, so we cannot determine the value of
$\epsilon_{i,j,i',j'} = [c_i,c_j]\cdot [c_{i'},c_{j'}]$ aside from
observing that it lies in the set $\{1,0,-1\}.$

 Now, the curve $u_1^{r_1}\cdots
u_{2g-2}^{r_{2g-2}}[c_{i'},c_{j'}]$ cannot possibly meet $[c_i,c_j]$
unless
$r_k\in\{\delta_{i,k},\delta_{j,k},-\delta_{i',k},-\delta_{j',k}\}$
for all $k.$ On the other 
hand, if the two curves do meet, then by symmetry, their intersection
numbers are determined by $[c_i,c_j]\cdot [c_{i'},c_{j'}]$, as 
indicated in the statement of the proposition.

The verification of equation~(\ref{eq:calculationparallel}) proceeds
in much the same way.  The only subtlety comes in checking the cases
$r_i = 0,$ in which the curves in question actually
have an entire segment in common.  But one can perturb one of the
curves so that they only meet at one endpoint of the segment; the
computation then follows from the usual symmetry.
\end{proof}

\subsection{The representation}

It will be useful for us to consider pointed versions of $\T_g$ and
$\K_g$.  We work with respect to the basepoint $x = p(y)\in \sg.$
Denote by $\tgs$ the group of components of the group of
basepoint-preserving diffeomorphisms of $\sg$ which act trivially on 
$H_1(\Sigma_g;\Z)$.  ``Forgetting the basepoint'' clearly gives a 
surjective homomorphism 
$\tgs\rightarrow\tg.$  Denote by $\kgs$ 
the subgroup of $\tgs$ generated by twists about bounding
curves which avoid the basepoint.  Again, the operation of forgetting
the basepoint induces a surjection $\kgs\rightarrow\kg.$

Recall that $K=\pi_1(Y)$.  Note that 
since $K$ is not a characteristic subgroup of $\pi_1(\sg),$ an
arbitrary mapping class need not lift to $Y$.  In fact, there are even
elements of $\T_g$ which don't lift to $Y$.  However, we have the following.

\begin{proposition}
Each element of $\kgs$ has a 
lift to a basepoint-preserving diffeomorphism of $Y$ which is
unique up to basepoint-preserving isotopy.
\end{proposition}

\begin{proof}
The uniqueness is clear.  Moreover, by the universal lifting
property for covering maps, the collection of basepoint-preserving
mapping classes that admit such a lift constitutes a subgroup.  Thus,
we need only verify the result for Dehn twists about bounding curves, as
these generate $\kgs$.

To this end, let $C$ be a bounding curve on $\sg,$ and denote by $t_C$ the
twist about $C.$  Since $p$ is an abelian cover, $C$ lifts to a simple
closed curve in $Y$.  Consider the map $\widetilde{t}_C$, which is a
simultaneous Dehn twist about all the lifts of $C.$  This obviously
constitutes a lift of $t_C.$
\end{proof}

These observations are enough to give us our main tool.  Henceforth $C$ 
will denote an arbitrary bounding curve in $\sg,$ and
$\widetilde{C}$ will denote a lift of $C$ to $Y$.  The homology class of
$\widetilde{C}$ will be written
$$c+\sum
m_{p_2,\dots,p_g,q_2,\dots,q_g}s_2^{p_2}\cdots s_g^{p_g}t_2^{q_2}\cdots
t_g^{q_g} a_1 +
n_{p_2,\dots,p_g,q_2,\dots q_g}s_2^{p_2}\cdots
s_g^{p_g}t_2^{q_2}\cdots t_g^{q_g} b_1$$
where $c\in W$ (recall $W$ was defined in the statement of 
Proposition \ref{Pidenthom}), the sum is taken over all integers
$p_2,\dots,p_g,q_2,\dots,q_g,$ and 
the $m$'s and $n$'s are integral coefficients, all but finitely many
of which vanish.  To simplify the notation, we will use underlined symbols
to refer to $(g-1)$-tuples of objects indexed by the set
$\{2,\dots,g\}.$  For example, $\underline{p}$ will stand for
$p_2,\dots,p_g$, the symbol $\underline{s}$ will stand 
for $s_2,\dots,s_g$ and, crucially,
binary operations on underlined quantities will be performed
componentwise, so that
$\underline{s}^{\underline{p}} = s_2^{p_2}\cdots s_g^{p_g}.$

We are now ready to lift the action of $\K_{g,\ast}$.

\begin{proposition}
\label{Pfixes}
The operation which associates to an element of $\kgs,$ the action of
its lift to $Y$ on $H_1(Y;\Z)$ gives rise to a representation
$$\widetilde{\rho}:\kgs\rightarrow \Aut_{\lg}(H_1(Y;\Z))$$
\end{proposition}

\begin{proof}
We must check that $\widetilde{\rho}$ takes composition to
multiplication and that its image respects the $\lg$-action on
$H_1(Y;\Z)$.   The
former condition follows from the uniqueness up to isotopy 
of lifts; the latter holds because for any bounding curve $C$ in $\sg,$ the set
of all lifts of $C$ to $Y$ is $\lg$-invariant. 
\end{proof}

\subsection{Reducing dimension}

The representation $\widetilde{\rho}$ is quite complicated, because
$H_1(Y;\Z)$ is a rather large module.  We 
instead would like to work with a $2$-dimensional
$\lg$-representation.  We will achieve this by proving that
$\widetilde{\rho}$ 
contains a large subrepresentation, namely $W$, 
that we will be able to ignore.  In order to do this we first need to
analyze the image under $\widetilde{\rho}$ of a twist about a bounding curve.

\begin{proposition}
\label{Pdivides}
Let $C$ be a bounding curve on $\sg$ and $1\leq i< j\leq 2g-2.$  Then
$$\tilde{\rho}(t_C)([c_i,c_j]) = [c_i,c_j] + d$$ where $d$ can be
written as a sum of terms each of which is divisible by
$(u_{k} - 1)(u_{l} - 1)$ for some $1\leq k<l\leq 2g-2.$
\end{proposition}

\begin{proof}
We first assume that $\tilde{C} = [c_{i'},c_{j'}].$  
Recall that if
$\beta = \{\beta_k\}$ is a family of mutually disjoint and nonisotopic 
simple closed curves on a surface, and if 
$\alpha$ is another simple closed curve, then the homology
class of the twist $t_{\beta}(\alpha)$ of $\alpha$ about $\beta$ is 
\begin{equation}\label{eq:twist}
[t_{\beta}(\alpha)]=[\alpha] +
\sum_k(\alpha\cdot\beta_k)[\beta_k]\end{equation}

Now, if $\{i,j\} = \{i',j'\},$ then of course
$\tilde{\rho}(t_C)([c_i,c_j]) = [c_i,c_j].$  If, instead, $\{i,j\}\cap
\{i',j'\} = \emptyset,$ then
equation~(\ref{eq:calculation}) tells us that
$$\tilde{\rho}(t_C)([c_i,c_j]) = [c_i,c_j] + \epsilon_{i,j,i',j'}(u_i
- 1)(u_j - 1)(u_{i'}^{-1}-1)(u_{j'}^{-1} - 1)[c_{i'},c_{j'}]$$
which is of the desired form if we set $k=i$ and $l=j.$  Lastly,
suppose that $\{i,j\}\cap\{i',j'\}$ contains a single element, say
without loss of generality $i = i'.$  Then
equation~(\ref{eq:calculationparallel}) gives us
\begin{eqnarray*}
\tilde{\rho}(t_C)([c_i,c_j]) & =& [c_i,c_j] +
\epsilon_{i,j,,j'}\left(-u_i^{-1}+1-u_i\right) \left(u_j - 
1\right)\left(u_{j'}^{-1} - 1\right)[c_{i'},c_{j'}] \\
& =& [c_i,c_j] +
\epsilon_{i,j,,j'}u_{j'}^{-1}\left(u_i^{-1}-1+u_i\right) \left(u_j -  
1\right)\left(u_{j'} - 1\right)[c_{i'},c_{j'}]
\end{eqnarray*}
which again gives us what we want, with $k = j$ and $l = j'.$

The general case follows from this calculation by the linearity
present in 
equation~(\ref{eq:twist}) along with the vanishing of
equation~(\ref{eq:vanishing}).
\end{proof}

We now use Proposition \ref{Pdivides} 
to find a substantially smaller representation of
$\kgs.$  Denote by $\lt$ the Laurent series ring $\Z[t,t^{-1}],$ and 
define $\Phi:\lg\longrightarrow\lt$  
by $\Phi(s_i)=1$ for $2\leq i\leq g$ and 
$$
\Phi(t_i)=\left\{
\begin{array}{ll}
t & \mbox{if $i=2$}\\
1 &\mbox{if $3\leq i\leq g$}
\end{array} \right.
$$

The homomorphism $\Phi$ induces a homomorphism 
$$\hat{\Phi}:\Aut_{\lg}(M)\longrightarrow \Aut_{\lt}(M\otimes_{\lg}\lt)$$
for any $\lg$-module $M.$  Now define 
$$\hat{\rho}:\kgs\rightarrow
\Aut_{\lt}\left(H_1(Y;\Z)\otimes_{\lg}\lt\right)$$ by 
$$\hat{\rho}=\hat{\Phi}\circ\widetilde{\rho}$$  

Recall now that $W\subset H_1(Y;\Z)$ was defined in the statement
of Proposition~\ref{Pidenthom}.

\begin{corollary}
The representation $\hat{\rho}$ becomes trivial when restricted to
$W\otimes_{\lg}\lt.$
\end{corollary}

\begin{proof}
Proposition~\ref{Pdivides} guarantees that for any bounding curve
$C$ and for any $1\leq i<j\leq 2g-2,$ that 
$$\widetilde{\rho}(t_C)([c_i,c_j]) = [c_i,c_j] + d$$ 
where $d$ is a sum
of terms each of which is
divisible by $(u_k - 1)(u_l - 1)$ for some $1\leq k< l\leq 2g-2.$
Since $k\neq l,$ at least one of $u_k$ and 
$u_l$ is not equal to $t_2,$ so it must be the case that each of the
summands of $d$ vanishes
when we tensor with $\lt.$  Thus, $\hat{\rho}(t_C)([c_i,c_j]) =
[c_i,c_j].$  The desired result then follows from the fact that the
$t_C$ generate $\kgs$ and the $[c_i,c_j]$ generate $W\otimes_{\lg}\lt.$
\end{proof}

We are now able to define the representation that will actually
allow us to prove our result.  Since the image of
$\hat{\rho}:\kgs\rightarrow \Aut_\lt\left(H_1(Y;\Z)\otimes_{\lg} \lt\right)$
fixes
$W\otimes_{\lg}\lt,$  we may pass to a quotient representation
$$\check{\rho}:\kgs\longrightarrow
\Aut_{\lt}\left[\frac{H_1(Y;\Z)\otimes_{\lg}\lt}{W\otimes_{\lg}\lt}\right]
\approx \GL_2(\lt)$$
where the last isomorphism follows from Proposition~\ref{Pidenthom}.

The first thing we will need to know about $\check{\rho}$ is the
following.

\begin{proposition}
The image of $\check{\rho}$ is actually contained in $\SL(\lt)$ rather than
$\GL_2(\lt).$ 
 Moreover, for a bounding
curve $C$ on $\sg,$ we have 
\begin{equation}\label{eq:rho}\check{\rho}(t_C) = 
\Phi\left(\begin{array}{cc}
1+\sum
n_{\underline{i},\underline{j}} m_{\underline{p},\underline{q}}
\underline{s}^{\underline{p}-\underline{i}}
\underline{t}^{\underline{q}-\underline{j}} & 
- \sum m_{\underline{i},\underline{j}}
m_{\underline{p},\underline{q}} 
\underline{s}^{\underline{p}-\underline{i}}
\underline{t}^{\underline{q}-\underline{j}} \\ 
\sum n_{\underline{i},\underline{j}}
n_{\underline{p},\underline{q}} 
\underline{s}^{\underline{p}-\underline{i}}
\underline{t}^{\underline{q}-\underline{j}} &
1-\sum m_{\underline{i},\underline{j}}
n_{\underline{p},\underline{q}} 
\underline{s}^{\underline{p}-\underline{i}}
\underline{t}^{\underline{q} - \underline{j}}
\end{array}\right)\end{equation}
Furthermore, $\check{\rho}$ descends to a representation
$$\rho:\kg\longrightarrow\SL(\lt)$$
\end{proposition}

\begin{proof}
Observe that the statement that the image of $\check{\rho}$ lies
in $\SL(\lt)$ rather than $\GL_2(\lt)$ follows formally from
equation~(\ref{eq:rho}), so it suffices to verify that equality.  To
establish that, we compute before projecting to $\lt$ via $\Phi$ by
simply expanding out the summations 
$$(\widetilde{t}_C)_*(a_1) \equiv a_1 + \sum (a_1\cdot
\underline{s}^{\underline{i}}\underline{t}^{\underline{j}}\widetilde{C}) 
\underline{s}^{\underline{i}}\underline{t}^{\underline{j}}[\widetilde{C}]
\pmod{W}$$
and
$$(\widetilde{t}_C)_*(b_1) \equiv b_1 + \sum (b_1\cdot
\underline{s}^{\underline{i}}\underline{t}^{\underline{j}}\widetilde{C}) 
\underline{s}^{\underline{i}}\underline{t}^{\underline{j}}[\widetilde{C}]
\pmod{W}$$
using the formulas
$$a_1\cdot
\underline{s}^{\underline{i}}\underline{t}^{\underline{j}}\widetilde{C} =
n_{-\underline{i},-\underline{j}}$$
and
$$b_1\cdot
\underline{s}^{\underline{i}}\underline{t}^{\underline{j}}\widetilde{C} =
-m_{-\underline{i},-\underline{j}}$$

To verify the last
statement, consider an element $\eta$ of $\kgs$ that lies in the kernel of
the projection $\kgs\rightarrow\kg.$  Denote by $\widetilde{\eta}$ the
basepoint-preserving lift of $\eta$ to $Y.$  Since $\eta$ is
isotopic to the identity once we forget basepoints, $\widetilde{\eta}$
must be isotopic to a diffeomorphism covering the identity map on
$\sg.$  Thus, we must have an equation
$$\check{\rho}(\eta) = \left(
\begin{array}{cc}
\Phi\left(\underline{s}^{\underline{p}}\underline{t}^{\underline{q}}\right)
& 0\\ 
0 &
\Phi\left(\underline{s}^{\underline{p}}\underline{t}^{\underline{q}}\right)
\end{array}\right)$$
But in order for this to lie in $\SL,$ it must be the identity matrix,
so $\check{\rho}$ factors through the quotient $\kg$ of $\kgs$.
\end{proof}

\section{Amalgamated products and infinite generation}

Denote by $H$ the image of the homomorphism 
$\rho:\kg\to\SL(\lt)$.  Our goal is to prove that $H$ is not 
finitely generated.  We now describe how we will do this.

Consider the inclusion 
$\SL(\lt)\subset\SL(\Q[t,t^{-1}])$.  The field $\Q(t)$ obtained by
adjoining a free variable $t$ to the rational numbers is
equipped with a discrete valuation and contains $\Q[t,t^{-1}]$, so
one can apply the construction
of Bruhat-Tits-Serre to find a (locally infinite) 
simplicial tree on which $\SL(\Q[t,t^{-1}])$ acts by isometries.  The
Bass-Serre theory of groups acting on trees can then be applied (see
\cite{BM}, \S 5) to express $\SL(\Q[t,t^{-1}])$ as an amalgamated product:

\begin{equation}\label{eq:amalg}
\SL(\Q[t,t^{-1}]) \cong A *_{\U} B
\end{equation}

\noindent
where $A = \SL(\Q[t]),$ 
$$B = \left(\begin{array}{cc} 
t^{-1} & 0 \\
0 & 1
\end{array}\right)
A
\left(
\begin{array}{cc} 
t & 0 \\ 
0 & 1
\end{array}
\right)$$
and $\U = A\cap B$.

This decomposition allows one to apply the theory of graphs of groups to
obtain the following criterion, which is Proposition 5 in \cite{MM}.

\begin{proposition}[Criterion for infinite generation]
\label{Pnotfg}
Let $A *_{\U} B$ be an amalgamated product, and let $H$ be any 
subgroup.  Suppose there exist elements $M_k\in A\backslash\U$ and
$N_k\in B\backslash\U$ such that
\begin{enumerate}
\item $M_k N_k M_k^{-1} \in H;$ and
\item $(H\cap A)M_k \U \neq (H\cap A)M_l \U$ whenever $k\neq l$.
\end{enumerate}
Then $H$ is not finitely generated.
\end{proposition}

We apply Proposition \ref{Pnotfg} to the situation above, with 
$\SL(\Q[t,t^{-1}]) \cong A *_{\U} B$ and with 
$H=\rho(\kg)$.  Our goal now is to find matrices $M_k$ and $N_k$
satisfying the desired criterion.

\subsection{The elements $M_k$ and $N$}

For a positive integer $k$, we let 
$$M_k=
\left(\begin{array}{cc}
1 & 0 \\
k & 1
\end{array}\right)
\in \SL(\lt)$$
We also set 
$$N=\left(\begin{array}{cc}
1 & t - 2 + t^{-1} \\
0 & 1
\end{array}\right)
\in \SL(\lt)$$

We now verify that the first hypothesis of 
Proposition \ref{Pnotfg} holds in our case; here we are taking 
$N_k=N$ for all $k$.  

\begin{proposition}
\label{Pliesin}
For each $k\geq 1,$ the matrix $M_k N M_k^{-1}$ lies in $H.$
\end{proposition}

\begin{proof}
First of all, consider the simple closed bounding curve $C$ shown in 
Figure~1.  The figure is drawn so that the homology of the leftmost
handle of $\sg$ is spanned by $\{a_1,b_1\}.$
\begin{figure}[h!]\label{fig:curve1}
\centerline{
\epsfig{figure = 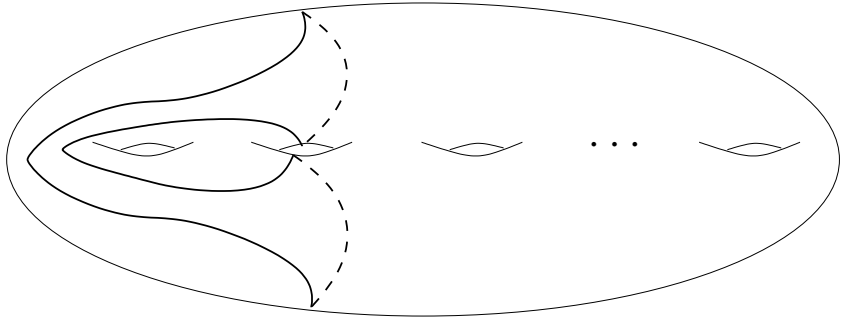}
}
\caption{The curve $C$}
\end{figure}
We now lift $C$ to $Y$; this is shown in Figure~2.
Here, each octogan with a handle coming out of it corresponds to a
single fundamental domain for the $\Z^{2g-2}$-action on $Y$; we have
drawn the two lifts of $C$ that meet the fundamental domain $X$.  In
general, of course, the base is a $4(g-1)$-gon; the figure
corresponds to the case $g=3$. 
\begin{figure}[h!]\label{fig:curve2}
\centerline{
\epsfig{figure = 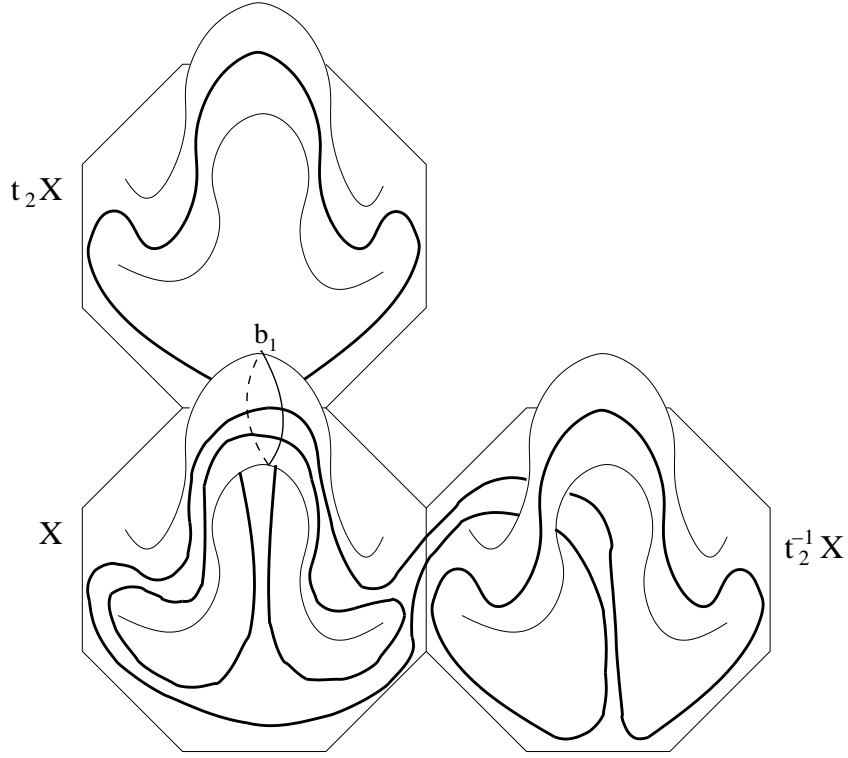}
}
\caption{The two lifts of $C$ that meet $b_1$.  }
\end{figure}
It is clear that no lifts of $C$ meet $a_1,$ so that $\rho(t_C)(a_1) =
a_1.$  Moreover, by twisting $b_1$ about the two curves shown in
Figure~2, one sees that $$\rho(t_C)(b_1) = b_1 +
\Phi(t_2 - 2 + t_2^{-1})a_1 = b_1 + (t - 2 + t^{-1})a_1$$ and therefore
that $\rho(t_C) = N.$

Secondly, since $b_1$ is in the kernel of the map
$\pi_1(\sg)\rightarrow\Z^{2g-2},$ the twist $t_{b_1}$ lifts to $Y$.  
Denoting by $T$ the simultaneous twist about all the lifts of $b_1$ to
$Y$, we see that $M_k = T_*^k.$  Set $C' = t^k_{b_1}(C)$.  We then have
\begin{eqnarray*}
M_k N M_k^{-1} &=& T_*^k \left(\widetilde{t}_C\right)_* T_*^{-k} \\
&=& \left(\widetilde{t}_{C'}\right)_*
\end{eqnarray*}
the last equality following from the general formula
$ft_af^{-1}=t_{f(a)}$, where $f$ is any mapping class and $t_a$ any
Dehn twist.  Since $C'$ bounds in $\sg,$ we see that
$M_k N M_k^{-1} = \rho(t_{C'}).$
\end{proof}

\subsection{Distinctness of double cosets}

The rest of this paper is devoted to proving the following.

\begin{proposition}\label{Pdoublecoset}
With the notation as above, we have
$$(H\cap A)M_k\U\neq (H\cap A)M_l\U$$ for all $k\neq l.$
\end{proposition}

Given Proposition \ref{Pdoublecoset}, whose proof we present in the next
section, we are now able to establish our main result, Theorem
\ref{theorem:main}.

\bigskip
\noindent
{\bf Proof of Theorem~\ref{theorem:main}. }We apply Proposition
\ref{Pnotfg} to the subgroup $H=\rho(\kg)$ of $\SL(\Q[t,t^{-1}]) \cong
A *_{\U} B$, 
with $M_k$ and $N_k=N$ as above.  First observe that 
$M_k\in A\backslash\U$ since $t$ does not divide $k,$
and $N\in B\backslash\U$ since $t-2+t^{-1}\not\in\Q[t].$
Therefore, in light of Propositions~\ref{Pliesin} and~\ref{Pdoublecoset},
Proposition~\ref{Pnotfg} implies that $H$ is not finitely generated.  As
$\kg$ surjects onto $H$, it is not finitely generated.
\endproof

Note that since $\kgs$ surjects onto $\kg$, it follows that $\kgs$ is
also not finitely generated.  

\section{The proof of Proposition \ref{Pdoublecoset}}

In this section we prove Proposition~\ref{Pdoublecoset}.  In order to do
this we will prove that 
the elements of $H=\rho(\kg)$ are of a very special form.  To state this
precisely, we will need the following. 

\begin{definition}[Balanced polynomials]
Let $f\in\Z[u_1^{\pm 1},\dots,u_n^{\pm 1}]$ be a Laurent polynomial in
$n$ variables over the integers.  We say that $f$ is {\em balanced} if
\begin{enumerate}
\item $f(1,1,\dots,1) = 0;$ and
\item for all $n$-tuples $(i_1,\dots,i_n)\in\Z^n,$ the coefficients of
$u_1^{i_1}\cdots u_n^{i_n}$ and $u_1^{-i_1}\cdots u_n^{-i_n}$ in $f$
are equal.

\end{enumerate}
\end{definition}

Parallel with a crucial observation of McCullough-Miller~\cite{MM}, we
have the following. 

\begin{proposition} \label{Pbalanced}
Each element of $H$ has the form
$$\left(\begin{array}{cc}
1 + P_1 & Q_1 \\
Q_2 & 1-P_2
\end{array}\right)$$
where $P_1,P_2,Q_1,$ and $Q_2$ are balanced.
\end{proposition}

\begin{proof}
Recall the map $\Phi:\lg\rightarrow\lt$ above.  
We begin by fixing an element $T = \rho(t_C)\in H$ and writing
$$T =
\left(\begin{array}{cc}
\Phi(1 + R_1) & \Phi(S_1) \\
\Phi(S_2) & \Phi(1 -R_2)
\end{array}\right)$$
where now the $R_i$ and $S_i$ lie in $\lg.$
Equation~(\ref{eq:rho}) gives us expressions for $R_1,R_2,S_1,$ and
$S_2$ in terms of the $m$ and $n$ coefficients.  Since the twist 
$t_C$ lies in the Torelli group $\T_g$, we have
$$\rho_*\left(\left(\widetilde{t}_C\right)_*\right)(a_1) = 
\rho_*(a_1) = a_1.$$  
But equation~(\ref{eq:rho}) tells us that 
$$\rho_*\left(\left(\widetilde{t}_C\right)_*\right)(a_1) = a_1 +
R_1(1,1,\dots,1)a_1 + S_2(1,1,\dots,1)b_1$$ 
and so $R_1(1,1,\dots,1)=S_2(1,1,\dots,1)=0$.  A
similar analysis of $b_1$ allows us to conclude that
$R_2(1,1,\dots,1)=S_1(1,1,\dots,1)=0.$ 

It follows via formal manipulations from equation~(\ref{eq:rho}) that
$S_1$ and $S_2$ 
also satisfy the other criterion for balancedness.  We now turn our
attention to $R_1$ and $R_2.$  Notice that for all $\underline{p}$ and
$\underline{q},$ we have 
\begin{eqnarray*}
0 & = &
\widetilde{C}\cdot\underline{s}^{-\underline{p}}\underline{t}^{-\underline{q}}
\\
& = &
\sum_{\underline{i},\underline{j}} m_{\underline{i},\underline{j}}
n_{\underline{i}+\underline{p},\underline{j}+\underline{q}} 
-
\sum_{\underline{i},\underline{j}} n_{\underline{i},\underline{j}}
m_{\underline{i}+\underline{p},\underline{j}+\underline{q}}
\end{eqnarray*}
From this, it follows that $R_1 = R_2,$ from which one can deduce
formally that $R_1$ is balanced.

Since it is clear that $\Phi$ takes balanced polynomials to balanced
polynomials, we have the desired property for elements of the form
$\rho(t_C).$  But the set of elements
of $\SL(\lt)$ of the desired form is evidently a subgroup, so the
result follows since the $t_C$ generate $\kg.$
\end{proof}

Following Lemma 7 in \cite{MM}, we will now see how 
Proposition~\ref{Pdoublecoset}  
follows rather formally from Proposition \ref{Pbalanced}.  

\begin{proof}[Proof of Proposition~\ref{Pdoublecoset}]
Suppose, that the $M_k$ and $M_l$ are in the same double coset, that
is, that we have a matrix equation
\begin{equation}\label{eq:matrices}
\left(\begin{array}{cc}
1 &0 \\
k & 1
\end{array}\right)
=
\left(\begin{array}{cc}
1+P_1 & Q_1\\
Q_2 &1-P_2
\end{array}\right)
\left(\begin{array}{cc}
1 & 0\\
l &1
\end{array}\right)
\left(\begin{array}{cc}
u & v\\
wt & z
\end{array}\right)
\end{equation}
with 
$$\left(\begin{array}{cc}
1+P_1 & Q_1\\
Q_2 & 1-P_2
\end{array}\right) \in H\cap A$$
and
$$\left(\begin{array}{cc}
u & v\\
wt & z
\end{array}\right) \in \U$$
By Proposition~\ref{Pbalanced},
we know that $P_1,$ $P_2,$ $Q_1,$ and $Q_2$ are balanced.  By the definition of
$A,$ they also lie in $\Q[t].$  Thus, they are constant and hence vanish.
Therefore, setting $t=0$ in equation~(\ref{eq:matrices}) gives
$$\left(\begin{array}{cc}
1 &0 \\
k & 1
\end{array}\right)
=
\left(\begin{array}{cc}
1 & 0\\
l &1
\end{array}\right)
\left(\begin{array}{cc}
u(0) & v(0)\\
0 & z(0)
\end{array}\right) =
\left(\begin{array}{cc}
u(0) & v(0)\\
lu(0) & lv(0) + z(0)
\end{array}\right)$$
which obviously implies that $k=l.$
\end{proof}

\medskip
\noindent
Dept. of Mathematics, University of Chicago\\
5734 University Ave.\\
Chicago, Il 60637\\
E-mail: {\texttt daniel@math.uchicago.edu}, {\texttt
farb@math.uchicago.edu} 

\end{document}